\documentclass{article}
\usepackage{amsthm}
\usepackage{amsfonts,amssymb,amsmath}
\usepackage{enumitem}
\usepackage{titlesec}
\usepackage{fouriernc}
\usepackage[T1]{fontenc}
\usepackage{chngcntr}
\counterwithin{figure}{section}

\newlist{gcases}{enumerate}{1}
\setlist[gcases,1]{
  label={{\it Case}~{\it \Alph*}.},
  topsep=0ex,
  leftmargin=0in,
  labelsep=.1in,
  itemindent=.7in,
  itemsep=0ex 
}
\newlist{tenumerate}{enumerate}{1}
\setlist[tenumerate,1]{
  label={(\arabic*)},
  topsep=0ex,
  leftmargin=.3in,
  labelsep=.1in,
  itemindent=0in,
  itemsep=0ex
}

\newlist{titemize}{enumerate}{1}
\setlist[titemize,1]{
  label={$\bullet$},
  topsep=0ex,
  leftmargin=.3in,
  labelsep=.1in,
  itemindent=0in,
  itemsep=0ex
}

\titleformat{\section}
  {\normalfont\bfseries}   
  {}                             
  {0pt}                          
  {Section \thesection\quad}    

\titleformat{name=\section,numberless}
  {\normalfont\bfseries}
  {}
  {0pt}
  {}

\newlength{\tabwidth}
\newlength{\tabheight}
\newlength{\tabrule}
\newlength{\tabwidthx}
\newlength{\tabheightx}

\def\gentabbox#1#2#3#4{\vbox to \tabheight{\setlength{\tabrule}{#3}%
  \setlength{\tabwidthx}{#1\tabwidth}\addtolength{\tabwidthx}{\tabrule}%

\setlength{\tabheightx}{#2\tabheight}\addtolength{\tabheightx}{-\tabheight}%
  \hbox to #1\tabwidth{%
 \hspace{-0.5\tabrule}\rule{\tabrule}{#2\tabheight}\hspace{-\tabrule}%
    \vbox to #2\tabheight{\hsize=\tabwidthx%
      \vspace{-0.5\tabrule}\hrule width\tabwidthx height\tabrule%
      \vspace{-0.5\tabrule}\vfil%
      \hbox to \tabwidthx{\hss#4\hss}%
        \vfil\vspace{-0.5\tabrule}%
      \hrule width\tabwidthx height\tabrule\vspace{-0.5\tabrule}}%
 \hspace{-\tabrule}\rule{\tabrule}{#2\tabheight}\hspace{-0.5\tabrule}}%
  \vspace{-\tabheightx}}}
\def\genblankbox#1#2{\vbox to \tabheight{\vfil\hbox to
#1\tabwidth{\hfil}}}
\def\tabbox#1#2#3{\gentabbox{#1}{#2}{0.4pt}{\strut #3}}

\catcode`\:=13 \catcode`\.=13 \catcode`\;=13
\catcode`\>=13 \catcode`\^=13
\def:#1\\{\hbox{$#1$}}
\def.#1{\tabbox{1}{1}{$#1$}}
\def>#1{\tabbox{2}{1}{$#1$}}
\def^#1{\tabbox{1}{2}{$#1$}}
\def;{\genblankbox{1}{1}\relax}
\catcode`\:=12 \catcode`\.=12 \catcode`\;=12
\catcode`\>=12 \catcode`\^=7

\newcommand\T{\mathbf{T}}
\newcommand\U{\mathbf{U}}

\begin{document}

\noindent
{\bf \large A family of operators generating domino tableaux of a fixed shape and a decomposition of left cells into isotypic components}

\vspace{.3in}
\noindent
W.~M.~MCGOVERN \\
\vspace{.1in}
{\it \small Department of Mathematics, University of Washington, Seattle, WA, 98195, USA
}

\section*{Introduction}
The purpose of this paper is to present and rework the results of \cite{H97}, exhibiting a set of operators defined on certain pairs of domino tableaux of the same shape and sending them to other such pairs with the same right tableau, such that given any two pairs of domino tableaux with the same right tableau there is a composition of operators sending the first pair to the second.  Using these operators we will also (finally) really correct the statements and proofs of Theorem 4.2 in \cite{M96} and Theorems 2.1 and 2.2 of \cite{M'96}, using these results to give explicit bases for the isotypic components of a Kazhdan-Lusztig left cell in types $B$ and $C$, as in \cite[Theorem 1]{M99}. (The statements of Lemma 1 and Theorem 1 of \cite{M99} are correct but the proofs are not; we will give correct proofs here.)  A later paper will treat the case of type $D$.

\section{Operators on tableau pairs}
We will use the notation of the three parts of Garfinkle's series of papers on the classification of primitive ideals in types $B$ and $C$.   We begin by defining the operators on domino tableaux that we will be using.   Given any pair $\alpha,\beta$ of adjacent simple roots of the same length, the wall-crossing operator $T_{\alpha\beta}$ is defined in \cite[2.1.10]{G92} on tableaux having $\beta$ but not $\alpha$ in their $\tau$-invariants.  The definition is extended to ordered pairs $(\T_1,\T_2)$ of domino tableaux of the same shape by decreeing that $T_{\alpha\beta}(\T_1,\T_2) = (T_{\alpha\beta}(\T_1),\T_2)$.  If instead $\alpha,\beta$ are adjacent but have different lengths, then there are two operators $T_{\alpha\beta}^L,T_{\alpha\beta}^R$ defined in \cite{G92}; they are defined on certain pairs of tableaux of the same shape.  We will use the operator $T_{\alpha\beta}^L$; if this operator is defined on a pair of tableaux, it sends this pair to one or two such pairs.  If $(\T_1,\T_2)$ is a tableau pair in the domain of $T_{\alpha\beta}^L$, then we define $U_{\alpha\beta}^L(\T_1,\T_2)$ to consist of the pair or pairs in $T_{\alpha\beta}(\T_1,\T_2)$ whose tableaux have the same shape as $\T_1$ and $\T_2$.  Thus every tableau pair in $U_{\alpha\beta}^L(\T_1,\T_2)$ has right tableau $\T_2$.  In type $C$, if the $1$- and $2$-dominos in $\T_1$ form a $2\times 2$ box, then the positions of both are replaced by their transposes to get the left tableau $\T_1'$ of one pair $(\T_1',\T_2)$ in $U_{\alpha\beta}^L(\T_1,\T_2)$.  If in addition the $2$-domino of $\T_1'$ lies in a closed cycle, then one moves through this closed cycle in $\T_1'$, leaving $\T_2$ unchanged, to produce a second pair $(\T_1'',\T_2)$ in $U_{\alpha\beta}^L(\T_1,\T_2)$.  Furthermore if $(\T_1',\T_2)$ is a pair in $U_{\alpha\beta}^L(\T_1,\T_2)$, then $(\T_1,\T_2)$ is a pair  in $U_{\beta\alpha}^L(\T_1',\T_2)$.  In type $B$ a similar recipe applies, except that if the $1$- and $2$-dominos in $\T_1$ form a hook (rather than a $2\times 2$ box), then each is replaced by its transpose to get one pair in $U_{\alpha\beta}^L(\T_1,\T_2)$; as before, if the $2$-domino now lies in a closed cycle, one moves through this closed cycle, leaving $\T_2$ alone, to produce a second pair in $U_{\alpha\beta}^L(\T_1,\T_2)$.  

As pointed out in \cite{H97}, compositions of operators $T_{\alpha\beta}$ and $U_{\alpha\beta}^L$ are not sufficient to generate all tableau pairs $(\T_1,\T_2)$ with a fixed right tableau $\T_2$.  One needs in addition four families of operators defined in an ad hoc way.  In type $C$, call the shape of a domino tableau a {\sl quasi-staircase} if it takes one of the forms $\sigma_n=(2n+1,\ldots,n+3,n+1,n+1,n-1,n-2,\ldots,1)$ or $\tau_n=(2n+2,\ldots,n+4,n+3,n+1,n+1,\linebreak n-1,n-2\ldots,1)$ for some $n\ge2$; thus a quasi-staircase shape is obtained from one of the two smallest such shapes $(5,3,3,1)$ and $(6,5,3,3,1)$ by adding hooks in the upper left corner so as to make both the first row and first column have length two more than they had before.   For every $n\ge2$, fix domino tableaux $\tilde\T_n,\tilde\T_n'$ of  respective shapes $\sigma_n,\tau_n$ such that the largest domino is vertical and is located at the end of the two rows of size $n+1$, while the next largest domino is horizontal and is located at the end of the row just above these two.  Define tableaux $\U_n,\U_n'$ similarly, interchanging the two largest dominos in $\tilde\T_n,\tilde\T_n'$, respectively.  Define operators $S_n$ (resp. $S_n'$) on tableaux pairs $(\T_1,\T_2)$ such that the first $n^2+n$ (resp. first $(n+1)^2$) dominos of $\T_1$ form a tableau $\tilde\T_n$ or $\U_n$ (resp. $\tilde\T_n'$ or $\U_n'$) by interchanging the two largest dominos in the copy of $\tilde\T_n$ or $\U_n$ (resp. $\tilde\T_n'$ or $\U_n'$) while leaving all other dominos in both $\T_1$ and $\T_2$ unchanged.  Similarly let ${}^tS_n, {}^tS_n'$ be the transposes of these operators, taking the transpose of any tableau pair in the respective domains of $S_n,S_n'$ to the transpose of its respective images under $S_n,S_n'$.  As in \cite{G93}, given a sequence $\Sigma$ of operators $T_{\alpha\beta},U_{\alpha\beta}^L,S_n,S_n',{}^tS_n$, and ${}^tS_n'$ there is a corresponding composition $T_\Sigma$ of operators taking a tableau pair $(\T_1,\T_2)$ in its domain to a set of tableau pairs $(\T_1',\T_2')$ such that $\T_2=\T_2'$ in all cases.

In type $B$ the new operators are defined similarly, except that the quasi-staircase shapes take the form $(2n+1,\ldots,n+3,n+2,n,n,n-2,n-3,\ldots,1)$ (and are tiled by $(n+1)^2$ dominos) or $(2n+2,\ldots,n+5,n+4,n+2,n+2,n,n-1,\ldots,1)$ (and are tiled by $n^2+2n$ dominos).  The operators $S_n,S_n',{}^tS_n,{}^tS_n'$ are then defined as in the previous paragraph, as is the composition $T_\Sigma$ of a sequence $\Sigma$ of operators.  Note that no shape of the form $(n,n-1,\ldots,k+2,k,k,k-1,\ldots,1)$ supports a domino tableau.

\section{Transitivity of the action on tableau pairs}

Following the arguments in \cite{H97}, we now prove the first of our two  results.

\newtheorem*{thm1}{Theorem 1}
\begin{thm1}
Given two pairs $(\T_1,\T_2),(\T_1',\T_2')$ of domino tableaux of the same shape such that $\T_2=\T_2'$, there is a sequence $\Sigma$ of operators $T_{\alpha\beta},U_{\alpha\beta}^L,S_n,S_n',{}^tS_n$, and ${}^tS_n'$ such that $(\T_1',\T_2')$ is one of the pairs in $T_\Sigma(\T_1,\T_2)$.
\end{thm1}

\begin{proof}
The proof is by induction on the number of dominos in $\T_1$, the  base case where this number is 1 being trivial.  Assume first that we are in type $C$.  Arguing as in the proof of \cite[3.2.2]{G93} we quickly reduce  the theorem to showing that given any tableau pair $(\T_1,\T_2)$ and an extremal position $P'$ in $\T_1$ there is a sequence $\Sigma$ of operators 
such that some pair $(\T_1',\T_2')$ in $T_\Sigma(\T_1,\T_2)$ has the largest domino of $\T_1$ in position $P'$.  To prove this we argue as in \cite {H97} and \cite[3.2.2]{G93}, breaking the argument into cases as was done there (and using the same notation).  Let $P=\{S_{i,j},S_{i,j+1}\}$ be the position of this domino in $\T_1$.   Cases A through G carry over without difficulty.  In Case H we must rework the cases $\rho_{i+1}(\T_1) = j-1$ and add the case $\rho_{i+1}(\T_1) = j-2$.  As in previous cases we find a position $P_1$ which is extremal in $\T_1\setminus \{P,P'\}$ (the tableau obtained from $\T_1$ by removing the positions $P$ and $P'$).  If $\rho_{i+2}(\T_1) = j-1$, then let $r=\kappa_{j-1}(\T_1)$ and set $P_1 = \{S_{r-1,j-1},S_{r,j-1}\}$.  If $p_{i+2}(\T_1) = j-2$, then there must be one more row or column as we have accounted for seven boxes so far.  If $i\ge2$, then let $u=\rho_{i-1}(\T_1)$ and set $P_1=\{S_{i-1,u-1},S_{i-1,u}\}$.  If $j\ge4$, then let $v=\kappa_{j-3}(\T_1)$ and set $P_1=\{S_{v-1,j-3},S_{v,j-3}\}$.  (If both conditions hold, then either choice of $P_1$ works.) We now prove the result as in case D.  The hardest case is case I; this is where the $S_n,S_n',{}^tS_n$, and ${}^tS_n'$ operators come into play.  Set $P'=\{S_{i-1,j+1},S_{i,j+1}\}$. 

 If $i=2$ and $j=1$, the $U_{\alpha\beta}^L$ operator suffices.  If $i=j=2$, set $p=\kappa_{j-1}(\T_1)$.    Then $p$ must be even.  If $p=2$, then an $F$-type interchange (defined in \cite[2.1.10]{G92} and implemented via a $T_{\alpha\beta}$ operator) does the job.  The case $p=4$ can be checked directly; here the tableaux have shape $(3,3,1,1)$.  The case $p\ge6$ is then handled by combining the transpose of case D with case H.  If $j\ge3$ there is no restriction on the parity of $p$.  If $p=2$ then let $P_1 = \{S_{1,j},S_{1,j+1}\},P_2=\{S_{1,j-1},S_{2,j-1}\}$ and argue as in subcase (a) of case I; here another $F$-type interchange arises.  If $p\ge5$, then let $P_1 = \{S_{p-1,j-1},S_{p,j-1}\}$ and proceed as above.  For $p=3$ or 4 let $q=\kappa_{j-2}(\T_1)$.  If $q=p$, then set $P_1=\{S_{p,j-2},S_{p,j-1}\}$ and apply case D, and if $q\ge p+2$ then argue as in the case $p\ge i+3$, using case G instead of H.  This leaves the cases $p=3$ or 4 and $q=p+1$.  Then $j\ne3$.  If $\kappa_{j-3}(\T_1) = \kappa_{j-2}(\T_1)$ or $\kappa_{j-3}(\T_1)\ge \kappa_{j-2}(\T_1)+2$, then argue as above.  Otherwise we must have $j\ne4$.  We continue in this fashion, eventually finding a column with an extremal position, since the overall tableau shape must be tilable by dominos.  
 
 If $i\ge3$ and $j=1$, then taking transposes in the argument for $i=2$ yields the desired result.  
 If $i\ge3$ and $j\ge2$, then let $r=\rho_{i-2}(\T_1)$.  If $p=3$ or $p\ge6$, then we argue as before, locating a column with an extremal position.  Transposes of above arguments treat the cases $r=j+1$ and $r\ge j+4$.  Thus we are reduced to the cases where $r=j+2$ or $j+3, p=4$ or $5, q = p+1$, and $\kappa_{j-3}(\T_1) = q+1$.  If $i=3, r= j+2$, then the overall tableau shape must be a quasi-staircase or the transpose of a quasi-staircase and the additional operators $S_n,S_n',{}^tS_n,{}^tS_n'$ do the job.  The remaining case J in \cite{G93} can be handled as in that paper.  This finishes all cases in type $C$.  The argument for type $B$ is analogous.
\end{proof}

We now \lq\lq enlarge" the operators $S_n,S_n'{}^tS_n$, and ${}^tS_n'$ to operators $T_n,T_n',{}^tT_n,{}^tT_n'$ that can take either one or two values.  Suppose we are given a tableau pair $(\T_1,\T_2)$ that can be moved through extended open cycles of $\T_1$ relative to $\T_2$ to produce a pair $(\T_1',\T_2')$ on which one of $S_n,S_n',{}^tS_n$, and ${}^tS_n'$, say $X$, is defined, and denote the enlarged operator by $X'$.  Let $d_1,d_2$ be the two largest dominos in $\T_1$ lying within the quasi-staircase or transposed quasi-staircase subshape in it, with $d_1>d_2$.  If $X=S_n'$ or ${}^tS_n'$, then let $\T_1'$ be obtained from $\T_1$ by interchanging dominos $d_1$ and $d_2$; the unique value of $X'(\T_1,\T_2)$ is $(\T_1'\T_2)$.  If $X=S_n$ or ${}^tS_n$ and dominos $d_1,d_2$ form a $2\times 2$ box in $\T_1$, then transpose the positions they occupy within the box to get a new tableau $\T_1'$.  Then one of the images of $X'(\T_1,\T_2)$ is $(\T_1',\T_2)$.  If the extended open cycle $e$ of $d_1$ in $
\T_1'$ relative to $\T_2$ contains $d_2$, then $(\T_1',\T_2)$ is the only image of $(\T_1,\T_2)$ under $X'$; if $e$ does not contain $d_2$, then move $(\T_1',\T_2)$ through $e$ to produce a second image $(\T_1'',\T_2')$ of $(\T_1,\T_2)$ under $X'$.  If $d_1,d_2$ occupy the positions specified by the definition of $X$, then set $(\T_1',\T_2)=X(\T_1,\T_2)$.  Then, as before, if the extended open cycle $e$ of $d_1$ in $\T_1'$ relative to $\T_2$ contains $d_1$, then we take $(\T_1',\T_2)$ to be the only image of $(\T_1,\T_2)$ under $X'$; otherwise, we move $(\T_1',\T_2)$ through $e$ to produce a second image $(\T_1'',\T_2')$ of $(\T_1,\T_2)$ under $X'$.  Finally, if dominos $d_1,d_2$ in $\T_1$ do not form a $2\times 2$ box or occupy the positions specified by the definition of $X$, then move $(\T_1,\T_2)$ through the extended open cycle of $d_2$ in $\T_1$ to produce a new pair $(\T_1',\T_2')$.  Then either transpose the positions of $d_1,d_2$ within the $2\times 2$ box they occupy in $\T_1$, or else apply $X$ to $(\T_1',\T_2')$, to get a new pair $(\T_1'',\T_2')$.  Take this pair to be the unique value of $(\T_1,\T_2)$ under $X'$. We will use these new operators in the next section to exhibit a number of maps from one left cell to another which intertwine the left action of the Weyl group.

\section{Decomposition of left cells}
We now recall the rule stated before Lemma 1 of \cite{M99} for constructing bases of isotypic components of Kazhdan-Lusztig left cells in Weyl groups $W$ of type $B$ or $C$ in terms of Kazhdan-Lusztig basis vectors $C_w$ of $W$.  Fix left cells $\mathcal C,\mathcal R$ of $W$ lying in the same double cell $\mathcal D$.  Let $x$ be the unique element of $\mathcal C\cap\mathcal R$ whose left tableau $T_L(x)$ has special shape.  Let $\sigma$ be the partition of $2n$ or $2n+1$ corresponding to a representation $\pi$ of $W$ occurring in both $\mathcal C$ and $\mathcal R$ and let $e_1,\ldots,e_r$ be the extended open cycles of $T_L(x)$ relative to $T_R(x)$ such that moving $T_L(x)$ through these open cycles produces a tableau of shape $\sigma$.  Given any $w\in\mathcal C\cap\mathcal R$, let $T_L(w)$ be obtained from $T_L(x)$ by moving through the extended open cycles  $f_1,\ldots,f_s$ (relative to $T_R(x)$).  Put $\sigma_w=\pm1$ according as an even or odd number of $f_i$ appear among the $e_j$.  Set $R_\sigma = \sum_{w\in\mathcal C\cap\mathcal R} \sigma_wC_w$.

\newtheorem*{thm2}{Theorem 2}
\begin{thm2}
The right or left $W$-submodule generated by $R_\sigma$ is irreducible and $W$ acts on it by $\pi$.
\end{thm2}

\begin{proof}
Since it is well known that both $\mathcal C$ and $\mathcal R$ (span vector spaces that) are multiplicity-free as $W$-modules, there must be some combination $R_\sigma'$ of the $C_w$ such that the right or left $W$-submodule generated by $R_\sigma'$ is isomorphic to $\sigma$, which is unique if one decrees that the coefficient of $C_x$ in $R_\sigma'$ is 1; moreover,  known facts about the asymptotic Hecke algebra of $W$ in the classical case imply that all coefficients of $R_\sigma'$ are 1 if $\pi$ is special while exactly half of them are 1 and the other half are $-1$ if $\pi$ is not special (see \cite{L87} and the proof of Theorem 1 in \cite{M99}).  More precisely, the coefficients of the $R_\sigma'$ as $\sigma$ runs over the partitions corresponding to representations appearing in both $\mathcal C$ and $\mathcal R$ are the rows of the character table of an elementary abelian 2-group.  It follows in particular that $R_\sigma' = R_\sigma$ for any cell intersection $\mathcal C\cap\mathcal R$ of size at most two.   The operators $T_{\alpha\beta}$ (for $\alpha,\beta$ of the same length) and $T^L_{\alpha\beta}$ (for $\alpha,\beta$ of different lengths) are well known to extend uniquely to $W$-equivariant linear maps from a left cell $\mathcal C'$ (regarded as a $W$-module) on which they are defined to its image $\mathcal C''$; more precisely, they are implemented by right multiplication by a suitable element in the Hecke algebra followed by projection to the relevant left cell.  (Whenever one of these operators sends a single tableau pair to two other such pairs, their linear extensions send a single Kazhdan-Lusztig basis vector to the sum of two such vectors.)  By \cite[2.3.4]{G92} they are also compatible with the formula for $R_\sigma$ given above in the sense that if a particular $R_\sigma$ coincides with $R_\sigma'$ and so transforms by $\pi$, then the same will be true of the image of $R_\sigma$ under any composition $\Sigma$ of maps $T_{\alpha\beta}$ or $T_{\alpha\beta}^L$ if $\Sigma$ is defined and nonzero on $R_\sigma$.  We similarly extend the operators $T_n,T_n',{}^tT_n$, and ${}^tT_n'$ to linear maps from any left cell on which they are defined to another such cell.

Now assume that $W$ is of type $C_6$ and consider a cell intersection $I=\mathcal C\cap\mathcal C^{-1}$ where $\mathcal C$ is a left cell represented by an element with left tableau $\tilde\T_1$ as chosen above for the shape $(5,3,3,1$, so that the irreducible constituents of $\mathcal C$ as a $W$-module are indexed by the partitions $(4,4,2,2),(4,3,3,2),(5,3,3,1)$, and $(5,4,2,1)$ of 12.  Denote by $x_p$ the unique element in the intersection $I$ whose left and right tableaux have shape $p$.  Then compositions of operators $T_{\alpha\beta}$ and $U_{\alpha\beta}^L$ act transitively on tableau pairs with a fixed right tableau not of shape $(5,3,3,1)$ or its conjugate $(4,3,3,1,1)$, so the argument in the proof of \cite[Theorem 1]{M99} applies to show that

\begin{align*}
R_{(4,4,2,2)}' =R_{(4,4,2,2)}= x_{(4,4,2,2)}+x_{(5,3,3,1)}+x_{(4,3,3,2)}+x_{(5,4,2,1)}\\ R_{(4,3,3,2)}' =R_{(4,3,3,2)}= x_{(4,4,2,2)}-x_{(5,3,3,1)}-x_{(4,3,3,2)} + x_{(5,4,2,1)}
\end{align*}
\vskip .2in
\noindent while either

\begin{align*}
R_{(5,3,3,1)}' = x_{(4,4,2,2)}+x_{(5,3,3,1)}-x_{(4,3,3,2)}-x_{(5,4,2,1)}\\
 R_{(5,4,2,1)}' = x_{(4,4,2,2)}-x_{(5,3,3,1)}+x_{(4,3,3,2)}-x_{(5,4,2,1)}
 \end{align*}
\vskip .2in
\noindent or else 

\begin{align*}
R_{(5,3,3,1)}' = x_{(4,4,2,2)}-x_{(5,3,3,1)}+x_{(4,3,3,2)}-x_{(5,4,2,1)}\\
R_{(5,4,2,1)}' = x_{(4,4,2,2)}+x_{(5,3,3,1)}-x_{(4,3,3,2)}-x_{(5,4,2,1)}
\end{align*}
\vskip .2in
\noindent But now if we take the basis elements for the Weyl group $W'$ of type $C_4$ corresponding to the tableau pairs consisting of the first four dominos of every tableau in all the pairs corresponding to elements of $I$ and label the resulting elements $y_q$ in type $C_4$ by partitions $q$ of 8 as we did the elements of $I$ by partitions of 12, we find that $y_{(2,2,2,2)} - y_{(3,2,2,1)}$ transforms by the representation corresponding to $(3,2,2,1)$ of $W'$, whose truncated induction to $W$ is the direct sum of the representations corresponding to $(5,4,2,1)$ and $(5,3,3,1)$.  In order to make $x_{(5,4,2,1)} - x_{(4,4,2,2)}$ transform by representations lying in this last truncated induced representation, we must have $R_{(5,3,3,1)}' = R_{(5,3,3,1)} = x_{(4,4,2,2)} +x_{(5,3,3,1)} - x_{(4,3,3,2)} - x_{(5,4,2,1)}$ and similarly for $R_{(5,4,2,1)}' = R_{(5,4,2,1)}$, as desired.  It follows that the operators $T_2,T_2',{}^tT_2$, and ${}^tT_2'$, extended to linear maps between left cells regarded as $W$-modules, are indeed equivariant for the left $W$-action.  Like the maps $T_{\alpha\beta}$ and $T_{\alpha\beta}^L$, they are compatible with the formula for $R_\sigma$.  Similar arguments show that the linear extensions of the other maps $T_n,T_n',{}^tT_n$, and ${}^tT_n'$ are also $W$-equivariant and compatible with the formula for $R_\sigma$.  Now we have enough $W$-equivariant maps between left cells to validate the proof of \cite[Theorem 1]{M99}.  A parallel argument handles Weyl groups of type $B$.
\end{proof} 

Now we can correct the statements of Theorem 4.2 in \cite{M96} and Theorems 2.1 and 2.2 in \cite{M'96}.  All three of these theorems are corrected and superseded by the following result:  {\sl given any two left cells $\mathcal C_1,\mathcal C_2$ in a Weyl group $W$ of type $B$ or $C$ that have a representation $\pi$ of $W$ in common, there is a composition $\Sigma$ of maps $T_{\alpha\beta},T_{\alpha\beta}^L,T_n,T_n',{}^tT_n,{}^tT_n'$ from $\mathcal C_1$ to $\mathcal C_2$ whose restriction to the copy of $\pi$ in $\mathcal C_1$ maps it isomorphically onto the corresponding copy of $\pi$ in $\mathcal C_2$}.  This follows since there are elements $w_1,w_2$ of $\mathcal C_1,\mathcal C_2$, respectively, such that the left tableaux of $w_1$ and $w_2$ both have shape the partition corresponding to $\pi$ and lying in the same right cell.  By Theorem 1 above, there is a composition $\Sigma$ of operators $T_{\alpha\beta},U_{\alpha\beta}^L,S_n,S_n',{}^tS_n$, and ${}^tS_n'$ mapping $w_1$ to $w_2$; the corresponding composition of the linear maps $T_{\alpha\beta},T_{\alpha\beta}^L,T_n,T_n',{}^tT_n$, and ${}^tT_n'$ does the trick.

\end{document}